\pdfoutput=1

\documentclass[12pt,draftcls,onecolumn]{IEEEtran}

\usepackage{graphics}         	 	% for pdf, bitmapped graphics files
\usepackage{epsfig}            		% for postscript graphics files
\usepackage{times} 				% assumes new font selection scheme installed
\usepackage[latin1]{inputenc}  		% latin fonts
\usepackage{amsmath}           		% assumes amsmath package installed
\usepackage{amssymb}           	% mathematic symbols
\usepackage{wasysym}           		% some extra characters
\usepackage{cite}				% smart citations
\usepackage{color} 				% colored text
\newcommand*\fvec[1]{\ensuremath{\mathbf{#1}}}                                	% fvec{} creates a bold symbol
\newcommand*\dx[3]{\frac{\partial^{#3} {#1}}{\partial {#2}^{#3}}}			% partial derivative, input: \dx{a}{b}{c}
\newenvironment{ea}[0]{\begin{eqnarray}} {\end{eqnarray}}                     	% just to shorten the input eqnarray
\newenvironment{ea*}[0]{\begin{eqnarray*}} {\end{eqnarray*}}                  	% just to shorten the input eqnarray*
\newenvironment{m}[0]{\begin{bmatrix}} {\end{bmatrix}}                        	% creates matrix with rectangular braces
\newcommand{\norm}[1]{\left\lVert{#1}\right\rVert}                     		% \norm{x} -->  ||x|| with proper spacing
\newcommand{\vectornorm}[1]{\left\lVert{#1}\right\rVert}                     	% \vectornorm{x} is old version of \norm{x}
\newcommand*\mc[0]{\mathcal}                                                  		% abbreviation for \mathcal
\newcommand*\mbb[0]{\mathbb}                                                  		% abbreviation for \mathbb
\newcommand*\mb[1]{\mbox{#1}}                                                 		% abbreviation for \mbox{}
\newcommand*\fao[1]{\left(\;\! \forall \: {#1} \right)}                       		% ( \forall [argument] ) - operator
\newcommand*\exo[1]{\left( \exists \, {#1} \right)}                           		% ( \exists [argument] ) - operator
\newcommand*\spano[1]{ \mbox{\rm span}\left\{ {#1} \right\} }            	% span{} operator with span in normal font
\newtheorem{theo}{Theorem}[section]       	% theorem enviroment theorem
\newtheorem{lem}{Lemma}[section]          	% theorem enviroment lemma
\newtheorem{cor}{Corollary}[section]      		% theorem enviroment corollary
\newtheorem{rem}{Remark}[section]         	% remark enviroment

\title{Geometric Analysis of the Formation Problem for Autonomous Robots}

\author{Florian D\"orfler and Bruce Francis% <-this % stops a space
\thanks{This work was supported by the Landesstiftung Baden-W\"urttemberg.}% <-this % stops a space
\thanks{Florian D\"orfler is with the Center for Control, Dynamical Systems and Computation, University of California at Santa Barbara, Santa Barbara, CA 93106, USA, {\tt\small dorfler@engineering.ucb.edu}. Bruce Francis is with the Electrical and Computer Engineering Department, University of Toronto, M5S3G4 Toronto, Canada, {\tt\small bruce.francis@utoronto.ca}}.
}

\begin{document}

%%%%%%%%%%%%%%%%%%%%%%%%%%%%%%%%%%%%%%%%%%%%%%%%%%%%%%%%%%%%%%%%%%%%%%%%%

% title page settings

\maketitle
\thispagestyle{empty}
\pagestyle{empty}

%%%%%%%%%%%%%%%%%%%%%%%%%%%%%%%%%%%%%%%%%%%%%%%%%%%%%%%%%%%%%%%%%%%%%%%%%

% -----------------------------------------------------------------------------------------------------%
%      Abstract
% -----------------------------------------------------------------------------------------------------%

\begin{abstract}
In the formation control problem for autonomous robots a distributed control law steers the robots to the desired target formation. A local stability result of the target formation can be derived by methods of linearization and center manifold theory or via a Lyapunov-based approach. It is well known that there are various other undesired invariant sets of the robots' closed-loop dynamics. This paper addresses a global stability analysis by a differential geometric approach considering invariant manifolds and their local stability properties. The theoretical results are then applied to the well-known example of a cyclic triangular formation and result in instability of all invariant sets other than the target formation.
\end{abstract}

%%%%%%%%%%%%%%%%%%%%%%%%%%%%%%%%%%%%%%%%%%%%%%%%%%%%%%%%%%%%%%%%%%%%%%%%%

% -----------------------------------------------------------------------------------------------------%
%      Keywords
% -----------------------------------------------------------------------------------------------------%

%% Note that keywords are not normally used for peerreview papers.
%\begin{IEEEkeywords}
%formation control, global stability analysis, hyperbolic manifolds
%\end{IEEEkeywords}

%%%%%%%%%%%%%%%%%%%%%%%%%%%%%%%%%%%%%%%%%%%%%%%%%%%%%%%%%%%%%%%%%%%%%%%%%

% -----------------------------------------------------------------------------------------------------%
%      Introduction
% -----------------------------------------------------------------------------------------------------%

\section{Introduction}\label{Section: Introduction}

\IEEEPARstart{T}{he} formation control of a network of autonomous mobile robots is an interesting instance of distributed control and motion coordination. In this setup  the autonomous robots have to be stabilized to a formation while each robot has only locally sensed information about the others. %Besides the imitation of biological flocking phenomena technological applications include mobile sensor networks forming antenna arrays, teams of UAVs performing surveillance missions, satellite formations for high-resolution imaging, and submarine swarms for oceanic exploration.

In the formation control problem graph theory plays a natural role, both to define  a formation and to describe the sensor relationships--who can ``see'' whom. Early work used the graph-theoretic concept of rigidity to construct undirected graphs \cite{Eren2002,OlfatiSaber02} suited for formation control. These concepts have been extended to directed graphs in \cite{Hendrickx}. An excellent reference reviewing the application of rigidity theory in formation control is \cite{Anderson09}. Recently  rigidity was employed as an analysis tool to show the stability of the desired target formation which is specified as an infinitesimally rigid framework \cite{Krick08,DoerflerECC09}. %,yu2009control,Summers2009}. 
Typically, a potential function approach is used to design distributed control laws, an approach that originally emerged for undirected graphs \cite{OlfatiSaber02} but has recently been extended to directed topologies \cite{Krick08,DoerflerECC09}. 
%An alternative formation control law steers the robot directly to a target point minimizing the distance constraints \cite{yu2009control,Summers2009}.% or the Jacobi overrelaxation iteration \cite{Cortes2009}.
In a potential function approach a natural Lyapunov function candidate is readily available and leads to an exponential stability result with a guaranteed region of attraction depending on the rigidity of the formation \cite{DoerflerECC09}. Local stability of the target formation can also be shown via methods of linearization and center manifold theory \cite{Krick08}, an approach that is also inherently related to rigidity. Neither of these approaches leads to global stability results since it is well known that there are various invariant sets of the robots' dynamics other than the target formation. A global stability analysis considering these sets has been carried out only for the benchmark example of a triangular formation \cite{Anderson,Cao1,Cao2,Cao2008,Smith} yielding convergence to the target formation from all but initially collinear formations. 

In the global stability analysis each of the references \cite{Cao1,Cao2,Anderson,Cao2008,Smith} follows a Lyapunov-type approach specific to the triangular formation which is not extendable to higher order formations. The present paper provides a tool independent of a Lyapunov function and based on differential geometry in order to rule out convergence of the robots to undesired equilibrium sets. These sets are parametrized as submanifolds embedded in the space of inter-agent positions, where the formation dynamics naturally evolve. A differential geometric stability tool for submanifolds is derived based on showing that the linearized vector field points away from these manifolds. This geometric result is based on purely algebraic computations and suffices to show instability of these submanifolds without guessing a Lyapunov function. In the application of this geometric method to the benchmark example of the triangular formation (with a cyclic sensor graph) we can confirm the results of \cite{Cao1,Cao2,Anderson,Cao2008,Smith}: initially non-collinear robots will be strictly bounded away from the set of collinear formations and converge exponentially to the desired target formation.

This paper is organized as follows: Section \ref{Section: Review of the Formation Control Problem} recalls the formation control problem for three robots. In Section \ref{Section: A Manifold Instability Theorem} the geometric method is derived and applied to the triangular formation in Section \ref{Section: Global Stability Analysis of the Triangular Formation} yielding a global stability result. Finally, some conclusions are drawn in Section \ref{Section: Conclusions}.

% -----------------------------------------------------------------------------------------------------%
%      Review of the Formation Control Problem for Three Robots
% -----------------------------------------------------------------------------------------------------%

\section{The Formation Control Problem for Three Robots}\label{Section: Review of the Formation Control Problem}

% -----------------------------------------------------------------------------------------------------%
\subsection{Review of the Setup}\label{Subsection: The Setup}

For our purposes an autonomous robot is a fully actuated vehicle in the plane that has no communication devices and is equipped only with an onboard camera. We assume that the robot's motion is modelled by the dynamics $\dot z_{i}= u_{i}$, where $z_{i} \in \mathbb{R}^{2}$ is the position of robot $i$ and $u_{i} \in \mathbb{R}^{2}$ is the control input. Altogether we consider three such robots and with the concatenated vectors\footnote{Vectors are written either as $n$-tuples or column vectors.} {$z = (z_{1} , z_{2} , z_{3})$} and $u = (u_{1} , u_{2} , u_{3})$ in  $\mathbb{R}^{6}$ the overall dynamics are $\dot{z}=u$.

\begin{figure}[b]
	\centering{
	\includegraphics[scale=0.38]{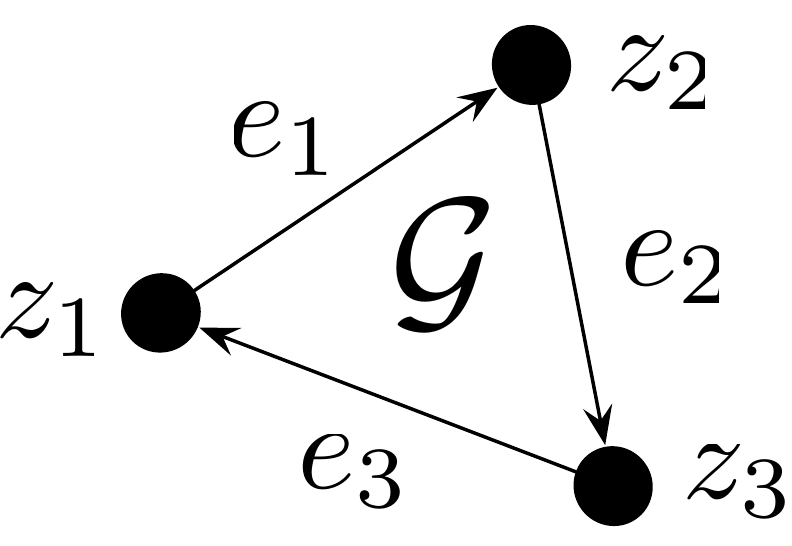}% draft version (12pt, single column, and double spacing)
	\caption{Illustration of the framework $(\mc G,z)$ together with the links $e$}
	\label{Fig: Framework of the triangle}
	}
\end{figure}

The sensing topology among the robots is specified by the cyclic \textit{sensor graph} $\mc G$, a directed graph with three nodes and three edges with clockwise orientation, as illustrated in Figure \ref{Fig: Framework of the triangle}. The nodes of $\mc G$ correspond to the robots, and we embed the graph into the plane as the {\it framework} $(\mc G,z)$. An edge $k$ from robot $i$ to robot $j$ corresponds to the \textit{link} $e_{k}=z_{j}-z_{i} \in \mbb R^{2}$ and means that robot $i$ can sense the relative distance and direction of  robot $j$ via its onboard camera. %Thus robot $i$ can sense the relative distance and direction of robot $j$.
 
We use the notation $e = (e_{1} , e_{2} , e_{3}) \in \mbb R^{6}$ for the concatenated vector of links, $I_{2}$ for the $2 \times 2$ identity matrix, and $0_{2}$ for the $2 \times 2$ zero matrix. With the block circulant {\it incidence matrix}
\begin{equation*}
	\hat H = 
	\left[ 	\begin{array}{rrr}
		-I_{2} & I_{2} & 0_{2} \\
		0_2 & -I_2 & I_2 \\
		I_2 & 0_2 & -I_2
	\end{array} \right] 
\end{equation*}
the links are obtained as $e=\hat H\,z$. The links are not independent, but subject to the constraint
\begin{equation}
	 e_{1}+e_{2}+e_{3} = \fvec 0
	 \label{eq: link space}
	 \,,
\end{equation}
where $\fvec 0 \in \mbb R^{2}$ is the vector of zero entries. The constraint \eqref{eq: link space} corresponds to the cycle $(1,2,3)$ in the graph $\mc G$ and defines a subspace in $\mbb R^{6}$ with normal vectors spanned by the columns of $(I_{2},I_{2},I_{2})$. We refer to this subspace as the \textit{link space} and denote it by $\mathrm{Im} \, \hat{H}$ (image of $\hat{H}$).

Given the sensor graph $\mc G$, a triangular formation is specified by a set of distance constraints $d_{k}>0$, $k \in \{1,2,3\}$, such that $ \| e_k \| = d_k$. Of course, the distance constraints have to be realizable, that is,  fulfill the triangle inequalities. The goal in formation control is to find a distributed control law $u_{i}=u_{i}(e_{i})$, that is, each control law can be implemented by onboard sensing, such that $z(t)$ converges as $t \rightarrow \infty$ and $\lim_{t} \| e_k(t) \| = d_k$ for all $k$. We refer to the set of all frameworks $(\mc G,z)$ fulfilling the distance constraints as the \textit{target formation}. 

In general, conditions to guarantee cohesion of the target formation and to stabilize the robots to it require a property called \textit{rigidity} of the target formation. Rigidity boils down to a rank condition on the \textit{rigidity matrix} $R_{\mc G}(e) = \mathrm{diag}\left(e_{i}\right)^{T} \hat H$: if the rank$R_{\mc G}(e)=3$ (it can't be more) then the {\it formation} $(\mc G,e)$ is said to be {\it infinitesimally rigid}. Infinitesimal rigidity is a generic property that holds in an open and dense set. In the triangular example all but collinear (and collocated) formations of robots are infinitesimally rigid and the additional necessary property of \textit{constraint consistence} \cite{Hendrickx} is also fulfilled. We do not further dwell on these properties but refer to \cite{Anderson09} reviewing rigidity theory and to \cite{Krick08,DoerflerECC09} %,yu2009control,Summers2009} 
relating it to sufficient stability conditions. 

Ideally the robots should converge to the target formation from any starting point. It is known that this goal cannot be achieved for every initial position $z(0)$, for example, the references \cite{Cao1,Cao2,Cao2008,Smith,Anderson,Krick08} show that three initially collinear robots cannot form a triangle. The objective of the present article is to provide a tool to find the exact region of attraction for the target formation.

% -----------------------------------------------------------------------------------------------------%
\subsection{A Potential Function Based Control Law}\label{Subsection: A Potential Function Based Control Law}

Typically a potential function approach is used to derive a distributed control law to tackle the formation control problem. For each robot a potential function is constructed that is zero whenever the robot has the desired distance from its neighbour and is positive when the distance constraints are violated. For robot $i$ define $W_{i} :\, \mbb R^{6} \to \mbb R$ as
$
W_{i}(z) = \frac{1}{4} \left( \norm{e_i}^2-d_{i}^2 \right)^{2}
$.
In order to minimize its potential, robot $i$ descends the gradient of the potential function, that is, $u_{i} = - \left[ \partial/\partial z_{i}\,W_{i}(z) \right]^{T}$. For notational convenience, we introduce the vector $\psi =(\psi_{1},\psi_{2},\psi_{3}) \in \mbb R^{3}$, where $\psi_{i} = \norm{e_i}^2-d_{i}^2$. The overall closed-loop $z$-dynamics are then
\begin{equation}
	\begin{bmatrix}
		\dot z_{1} \\ \dot z_{2}  \\  \dot z_{3}
	\end{bmatrix} 
	=
	\begin{bmatrix}
		e_{1}\,\left( \norm{e_1}^2-d_{1}^2 \right) \\
		e_{2}\,\left( \norm{e_2}^2-d_{2}^2 \right) \\
		e_{3}\,\left( \norm{e_3}^2-d_{3}^2 \right)
	\end{bmatrix} 
	=
	\begin{bmatrix}
		e_{1}\,\psi_{1} \\
		e_{2}\,\psi_{2} \\
		e_{3}\,\psi_{3}
	\end{bmatrix}
	\,,
	\quad z(0) = z_{0} \in \mbb R^{6} \,.
	\label{eq: z-dynamics} 
\end{equation}
Different approaches analyzing the $z$-dynamics in the state space $\mbb R^{6}$ have been proposed \cite{OlfatiSaber02,Krick08}. The target formation set in $\mathbb{R}^2$, i.e., the triangle, is invariant under rigid body motion. When lifted up to $\mathbb{R}^6$, the home of $z$, this set is non-compact. This complicates an analysis based on differential geometry, set stability or invariance concepts. In addition, the formation specification is in the link space. Fortunately, the target formation parametrized in the link space, 
\begin{equation*}
	\mc E_{e} = \{ e \in \mathrm{Im}\hat H:  \norm{e_{k}}=d_{k} , k=\{1,2,3\} \} \,,
\end{equation*}
is compact. For these obvious reasons we approach the stability analysis of the target formation in the link space. The closed-loop \textit{link dynamics}  resulting from the $z$-dynamics are
\begin{equation}
	\begin{bmatrix}
		\dot e_{1} \\ \dot e_{2}  \\  \dot e_{3}
	\end{bmatrix} 
	= 
	\begin{bmatrix}
		\dot z_{2}-\dot z_{1} \\ \dot z_{3}-\dot z_{2} \\ \dot z_{1}-\dot z_{3}
	\end{bmatrix} 
	=
	\begin{bmatrix}
		e_{2}\,\psi_{2} - e_{1}\,\psi_{1} \\
		e_{3}\,\psi_{3} - e_{2}\,\psi_{2} \\
		e_{1}\,\psi_{1} - e_{3}\,\psi_{3}
	\end{bmatrix}
	\,
	\quad e(0) = e_{0} = \hat H z_{0} \,.
	\label{eq: e-dynamics}
\end{equation}
The flow of the link dynamics in the link space $\mathrm{Im} \, \hat{H}$ will be denoted by $\Phi(t,e_{0})$.

% -----------------------------------------------------------------------------------------------------%
\subsection{A Preliminary Stability Result of the Target Formation}\label{Subsection: A Preliminary Stability Result of the Target Formation}

An intriguing approach to prove stability of $\mc E_{e}$ is to use the somewhat natural set-Lyapunov function candidate $V: \mathrm{Im} \, \hat{H} \to \mbb R$ defined as the sum of the potential functions
\begin{equation*}
	V(e) 
	=
	\sum\nolimits_{i=1}^{3} \frac{1}{4} \left( \norm{e_i}^2-d_{i}^2 \right)^{2} 
	=
	\frac{1}{4}\,\psi(e)^{T} \psi(e) \,.
\end{equation*}
The derivative of $V(e)$ along trajectories of the link dynamics can be compactly formulated as
\begin{equation}
	\dx{V}{e}{} \dot e 
	=
	- \norm{e_{1}\psi_{1}-e_{2}\psi_{2}}^{2} - \norm{e_{2}\psi_{2}-e_{3}\psi_{3}}^{2} - \norm{e_{3}\psi_{3}-e_{1}\psi_{1}}^{2}
	=
	- \psi^{T}\,R_{\mathcal{G}}(e)\,R_{\mathcal{G}}(e)^{T}\,\psi
	 \label{eq: dot V(e)} \,,
\end{equation}
%\begin{align}
%	\dx{V}{e}{} \dot e 
%	&=
%	- \norm{e_{1}\psi_{1}-e_{2}\psi_{2}}^{2} - \norm{e_{2}\psi_{2}-e_{3}\psi_{3}}^{2} - \norm{e_{3}\psi_{3}-e_{1}\psi_{1}}^{2}
%	\notag \\
%	&=
%	- \psi^{T}\,R_{\mathcal{G}}(e)\,R_{\mathcal{G}}(e)^{T}\,\psi
%	 \label{eq: dot V(e)} \,.
%\end{align}
where $R_{\mc G}(e)$ is the rigidity matrix. With the notation $\Omega(c) = \{e \in \mathrm{Im} \, \hat{H}  : \, V(e) \leq c \}$ for a sublevel set of $V(e)$ the following theorem can easily be derived from \eqref{eq: dot V(e)}:

\begin{theo}\label{Theorem: Semiglobal Exp. Stability}\cite[Theorem 5.1]{DoerflerECC09}
For every initial condition $e_{0} \in \mathrm{Im} \, \hat{H} $ the link dynamics \eqref{eq: e-dynamics} are forward complete and bounded in the compact sublevel set $\Omega(V(e_{0}))$, and their solution $\Phi(t,e_{0})$ converges to the largest invariant set contained in
\begin{equation*}
	\mc W_{e} 
	=
	\{ 
		e \in \Omega(V(e_{0})) :\, 
		\psi^{T}\,R_{\mathcal{G}}(e)\,R_{\mathcal{G}}(e)^{T}\,\psi = 0 
	\}
	\,.
\end{equation*}
Moreover, given $\rho>0$ such that for every $e \in \Omega(\rho)$ the formation $(\mc G,e)$ is infinitesimally rigid, for every initial condition $e_{0} \in \Omega(\rho)$ the set $\mathcal{E}_{e}$ is exponentially stable.% with respect to the link dynamics.
\end{theo}

By Theorem \ref{Theorem: Semiglobal Exp. Stability} the link dynamics converge either to the target formation $\mc E_{e}$ or the set $\mc W_{e} \setminus \mc E_{e}$, that is, the set of points in $\mc W_{e}$ where the matrix $R_{\mc G}(e)^{T}$ has a rank loss, spoken differently the set of non-rigid (i.e., collinear) formations. Locally the robots converge to the specified triangular formation with $\Omega(\rho)$ as guaranteed region of attraction. Note that $\Omega(\rho)$ is not necessarily a small set since rigidity is a generic property. As a result of the exponential convergence rate, the right-hand side of the $z$-dynamics \eqref{eq: z-dynamics} can be upper-bounded by exponentially decreasing signals and thus the positions also converge. Therefore, locally for every initial condition $z_0 \in \hat H^{-1}\left(\Omega(\rho)\right)$ the convergence of the robots to the formation is provable in  straightforward fashion \cite{DoerflerECC09}. 

Theorem \ref{Theorem: Semiglobal Exp. Stability} has a game-theoretic interpretation and also extends to a wider variety of graphs including undirected minimally rigid graphs \cite{DoerflerECC09}: for these graphs the only possible positive limit sets are the (locally stable) target formation and non-rigid formations. However, this result is only local and we are interested in the global behavior of the robots in the link space. Thus we have to find out the stability properties of the non-rigid sets. Such a global analysis for the triangular benchmark problem has been undertaken in \cite{Cao1,Cao2,Anderson} and for slightly different graphs in \cite{Smith,Cao2008} using problem-specific Lyapunov approaches. The next section provides a geometric method that allows an alternative approach by analyzing the linearized link dynamics only.

% -----------------------------------------------------------------------------------------------------%
%      A Manifold Instability Theorem
% -----------------------------------------------------------------------------------------------------%

\section{A Manifold Instability Theorem}\label{Section: A Manifold Instability Theorem}

The limit set $\mc W_{e}$ of the link dynamics  can be split into the target formation $\mc E_{e}$ and the set $\mc W_{e} \setminus \mc E_{e}$ of non-rigid limit sets. In order to show that $\mc W_{e} \setminus \mc E_{e}$ is not a positive limit set, it has to be shown that the vector field, the right-hand side of the link dynamics \eqref{eq: e-dynamics}, is pointing away from $\mc W_{e} \setminus \mc E_{e}$. This section formulates this idea in terms of differential geometry.

% -----------------------------------------------------------------------------------------------------%
\subsection{The Notion of Overflowing Invariance}\label{Subsection: The Notion of Overflowing Invariance}

Consider  the dynamical system
\begin{ea}
	\dot x = f(x) \;, \quad x(0)=x_{0} \in \mbb R^{n}
	\label{eq: Standard nonlinear system} \,,
\end{ea}where  $f: \, \mbb R^{n} \to \mbb R^{n}$ is a twice continuously differentiable vector field generating the flow $\Phi(t,x)$. In what follows,  $f_{x}(p)$ will denote the Jacobian of $f(x)$ at $x=p$. Let $\mc M$ be an $m$-dimensional differentiable submanifold $\mc M$ embedded in $\mbb R^{n}$ that is invariant w.r.t. \eqref{eq: Standard nonlinear system}, that is, for every $x_{0} \in \mc M$, $\Phi(t,x_{0}) \in \mc M$ for all $t \geq 0$. The normal and tangent space at $p \in \mc M$ are denoted as $N_{p}\mc M$ and $T_{p}\mc M$, and the normal and tangent bundles as $N\mc M$ and $T\mc M$. Geometrically speaking the invariance of $\mc M$ with respect to \eqref{eq: Standard nonlinear system} is equivalent to $f(p) \in T_{p}\mc M$ for all $p \in \mc M$.

The specification of $\mc M$ as an embedded submanifold allows us to identify a normal direction relative to $\mc M$. Given an $\epsilon > 0$, we can always construct a neighbourhood of $\mc M$ consisting of points $\tilde p \in \mbb R^n$ that are not further than $\epsilon$ away from $\mc M$ \cite[Theorem 6.17]{lee2003ism}. This can be seen as an embedding of the normal bundle $N\mc M$ into $\mbb R^{n}$ and we define the \textit{tubular $\epsilon$ neighbourhood}
\begin{equation*}
       \mc M_{\epsilon} 
       :=
	\{  
	\tilde p \in \mbb R^n : \, 
	\tilde p = p +  \bar\epsilon \, n_{p} ,  p \in \mathcal M  , n_{p} \in N_{p}\mc M , \norm{n_{p}} = 1 , 
	\bar \epsilon \in (0,\epsilon) 
        \} 
        \,.  
\end{equation*}
%\begin{align*}
%       \mc M_{\epsilon} 
%       :=
%	\{&  
%	\tilde p \in \mbb R^n : \, 
%	\tilde p = p +  \bar\epsilon \, n_{p} ,  p \in \mathcal M  , n_{p} \in N_{p}\mc M , \norm{n_{p}} = 1 , 
%	\notag \\
%	&\bar \epsilon \in (0,\epsilon) 
%        \} 
%     \,.  
%\end{align*}
We denote the boundary of the tubular $\epsilon$ neighbourhood $\mc M_{\epsilon}$ by $\partial \mathcal M_\epsilon$:
\begin{equation*}
    \partial \mathcal M_\epsilon 
    :=
    \{ 
    	\tilde p \in \mbb R^n : \, 
	\tilde p = p + \epsilon n_{p}  , p \in \mathcal M , n_{p} \in N_{p}\mc M , \norm{n_{p}} = 1
    \}  \,.
\end{equation*}
Let $\bar{\mc M}_\epsilon := \mc M \cup \mc M_\epsilon \cup \partial \mc M_\epsilon$ be the closure of $\mc M_{\epsilon}$.
Next we define the orientation of the vector field $f$ on $\partial \mc M_\epsilon$. Consider an $\epsilon > 0$, a point $p \in \mc M$, and a normal vector $n_{p} \in N_{p}\mc M$ of unit length. From this we construct the point $\tilde p \in \partial \mc M_{\epsilon}$ as $\tilde p = p + \epsilon \, n_{p}$. The inner product of the vector field $f(\tilde p)$ and the normal vector $n_{p}$ is then
\begin{equation}
    \left\langle f \left( \tilde p \right) ,  n_{p} \right\rangle  
    =
    \left\langle f \left( p + \epsilon n_{p} \right) ,  n_{p} \right\rangle 
    \label{eq: inner product}\,.
\end{equation}If the inner product \eqref{eq: inner product} is positive, then the vector field and the normal vector point in the same half space. We then say the vector field $f(\tilde p)$ is \textit{pointing strictly outward} at $\tilde p \in \partial \mc M_{\epsilon}$. Note that this property depends on $f$, $\epsilon$, $p$, and $n_{p}$. Consider a set $\Omega$ with $\mc M \cap \Omega \neq \emptyset$. If there exists an $\epsilon > 0$, such that for every $p \in \mc M \cap \Omega$ and for every $n_{p} \in N_{p}\mc M$ with $\norm{n_{p}}=1$ the vector field is pointing strictly outward, then we say $\mc M_{\epsilon}$ is \textit{overflowing invariant} in $\Omega$.

\begin{rem}
The term overflowing invariance is taken from \textit{Fenichel Theory}, which treats the stability properties of differentiable manifolds with boundaries \cite{Wiggins}. The invariant manifolds arising in our problem setup have no boundaries and thus this theory is not directly applicable.
\end{rem}

% -----------------------------------------------------------------------------------------------------%
\subsection{A Manifold Instability Result}\label{Subsection: A Manifold Instability Result}

The definition of overflowing invariance does not provide an easily checkable condition, since it depends on the, possibly nonlinear, vector field $f$ and the variables $\epsilon>0$, $p \in \mc M$, and $n_{p} \in$ $N_{p}\mc M$. Note that every embedded submanifold may be parameterized locally by the zero set of a smooth function \cite[Proposition 5.28]{lee2003ism}. In particular, consider the global case, where a continuously differentiable function $F:\, \mbb R^{n} \to \mbb R^{n-m}$ defines the zero set $\mc M := F^{-1}(\fvec 0)$. If rank$F_{x}(p)=n-m$ for all $p \in \mc M$, then $\mc M$ is an $m$-dimensional embedded submanifold, $F$ is said to be its {\it global defining function}, and the columns of the Jacobian $F_{x}(p)^{T}$ are a basis for $N_{p} \mc M$ \cite[Corollary 5.24, Lemma 5.29]{lee2003ism}. In this case, the idea to derive a checkable algebraic condition of overflowing invariance is to contract the tubular $\epsilon$ neighbourhood of $\mc M$ to a thin layer, in fact, to such a thin layer that the Taylor linearization of the vector field is valid. 

\begin{theo}\label{Theorem: Overflowing invariance}
Consider the vector field $f$ and an invariant embedded submanifold $\mc M := F^{-1}(\fvec 0)$ with the global defining function $F:\, \mbb R^{n} \to \mbb R^{n-m}$. Let $\Omega$ be a compact set with compact and non-empty intersection $\mc M \cap \Omega$, and consider for every $p \in \mc M \cap \Omega$ the matrix
\begin{equation}
    \Gamma(p) 
    =
     F_{x}(p)  \bigl( f_x(p) + f_x(p)^T \bigr) F_{x}(p)^{T} 
     \; \in \mathbb R^{(n-m) \times (n-m)} 
     \label{eq: overflowing invariance condition} \,.
\end{equation}
Assume that $\Gamma(p)$ is positive definite for every $p \in  \mc M \cap \Omega$. Then there exists $\epsilon^* > 0 $ such that, for every $\epsilon \in (0,\epsilon^{*}]$, the tubular $\epsilon$ neighbourhood $\mathcal M_{\epsilon}$ is overflowing invariant in $\Omega$.
\end{theo}

\begin{proof}
Let $\epsilon>0$ be arbitrary. We look at a point $\tilde p \in \partial \mc M_\epsilon$. By definition, it has the form
$
\tilde p = p + \epsilon n_{p}
$ 
for some $p \in \mc M$ and $n_{p} \in N_p \mc M$ with $\vectornorm{n_{p}}=1$. With $N_{p} \mc M = \mathrm{Im}F_{x}(p)^{T}$, $n_{p}$ can be parametrized as
$
n_{p} = F_{x}(p)^{T} \fvec c
$, 
where $\fvec c \in \mathbb{R}^{n-m}$. The inner product of $f\left(\tilde p\right)$ and $n_{p}$ is then % given by
\begin{equation*}
    \left\langle f \left(\tilde p \right) , n_{p} \right\rangle 
    =
    \left\langle f \left( p + \epsilon n_{p} \right)  ,  n_{p} \right\rangle
    =
    \left\langle f \left( p + \epsilon \, F_{x}(p)^{T} \fvec c \right) , F_{x}(p)^{T} \fvec c \right\rangle \,.
\end{equation*}
The ingredients $\mc M$, $\mc M_{\epsilon}$, and $f$ are illustrated in Figure \ref{Fig: Illustration of Overflowing Invariance} together with a trajectory.
\begin{figure}[b]
	\centering{
	\includegraphics[scale=0.34]{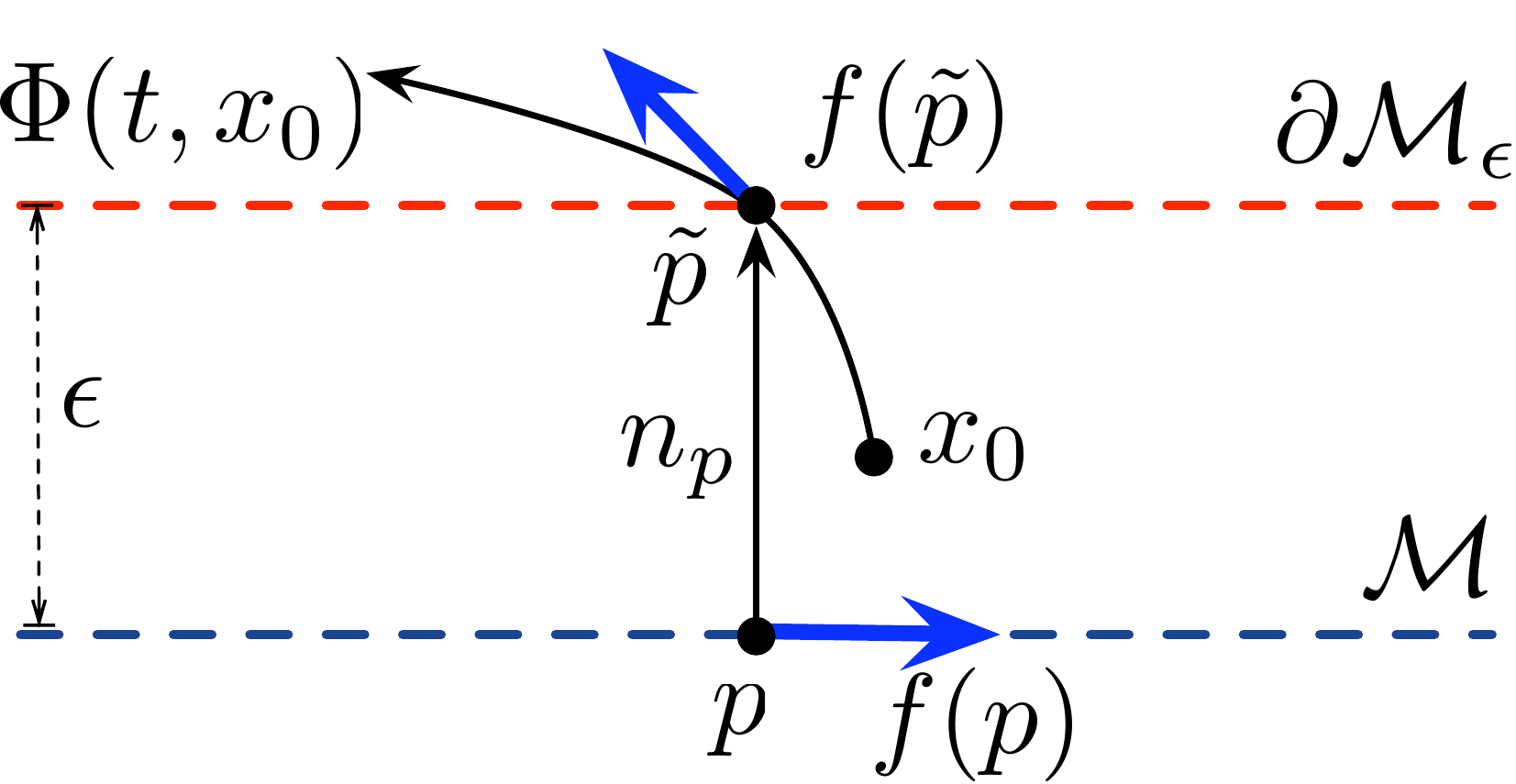} % draft version (12pt, single column, and double spacing)
	\caption{Qualitative illustration of $\mc M$, $\mc M_{\epsilon}$, $f$, and a trajectory $\Phi(t,x_{0})$ with $x_{0} \in \mc M_{\epsilon}$.}
	\label{Fig: Illustration of Overflowing Invariance}
	}
\end{figure}
We now expand $f \left( p + \epsilon \,F_{x}(p)^{T} \fvec c \right)$ in a Taylor series about $p \in \mc M$ and obtain for the inner product
\begin{equation}
    \left\langle f \left(\tilde p \right) , n_{p} \right\rangle 
    =
    \left\langle f(p)  ,  n_{p}  \right\rangle 
    + 
    \left\langle  \epsilon f_x(p)  n_{p} \,,\,  n_{p}  \right\rangle  
    + 
    \left\langle  R_3\left(p,\epsilon \right) , n_{p} \right\rangle 
    \label{eq: Taylor linearization of inner product 1} \,,
\end{equation}
where $R_3\left(p,\epsilon \right)$ is the Lagrange remainder of the Taylor series expansion and is of second order in $\epsilon$ %and $\norm{n(p)}$
 \cite[Theorem 4.1]{MultiVariableCalculus}. Note that the first term of \eqref{eq: Taylor linearization of inner product 1}  vanishes because $n_{p} \in N_{p}\mc M$ and $f(p) \in T_{p}\mc M$ due to invariance of the manifold $\mc M$. Thus equation \eqref{eq: Taylor linearization of inner product 1} simplifies to
\begin{equation}
    \left\langle f \left(\tilde p \right) , n_{p} \right\rangle 
    =
    \frac{\epsilon}{2} \, \fvec c^T \Gamma(p) \fvec c + \left\langle  R_3\left(p,\epsilon \right) , n_{p} \right\rangle  
    \label{eq: Taylor linearization of inner product 2} \,,
\end{equation}
where $\Gamma(p)$ is defined in \eqref{eq: overflowing invariance condition}. By definition $\mc M_{\epsilon}$ is overflowing invariant in $\Omega$ if the inner product \eqref{eq: Taylor linearization of inner product 2} is positive for every $p \in \mc M \cap \Omega$. 
If the  symmetric matrix $\Gamma(p)$ is positive definite, it is clear that we can obtain a positive inner product at every point $p \in \mc M \cap \Omega$ by choosing  $\epsilon$ sufficiently small at $p$. Let $\tilde \epsilon$ be such a sufficiently small $\epsilon$ at $p \in \mc M \cap \Omega$. Then we have
\begin{equation*}
    \fao{p \in \mc M \cap \Omega }
    \exo{\tilde \epsilon > 0}\;
    \frac{1}{2} \,\fvec c^T \Gamma(p) \fvec c 
    >
    \frac{1}{\tilde \epsilon} \left|\left\langle R_3\left(p,\tilde\epsilon \right) , n_{p} \right\rangle \right| 
    \,.
\end{equation*}
The right-hand side of the previous equation is upper bounded by the maximum Lagrange remainder, and by assumption, we have that $\Gamma(p)$ is positive definite for every $p \in \mc M \cap \Omega$:
\begin{equation}
    \fao{p \in \mc M \cap \Omega}
    \exo{\Gamma^{*} > 0 \,,\, R^*>0 }\;
    \Gamma(p) - \Gamma^*\,I_{n-m}
    \geq
    0
    \;,\;
    \tilde \epsilon R^* 
    \geq
    \frac{1}{\tilde \epsilon} \left|\left\langle  R_3\left(p,\tilde\epsilon \right) , n_{p} \right\rangle \right|  
    \label{eq: Bounds on Gamma and R} \,.
\end{equation}
%\begin{multline}
%    \fao{p \in \mc M \cap \Omega}
%    \exo{\Gamma^{*} > 0 \,,\, R^*>0 }\;
%    \\
%    \Gamma(p) 
%    \geq
%     \Gamma^*\,I_n
%     \;,\;
%    \tilde \epsilon \, R^* 
%    \geq
%    \frac{1}{\tilde \epsilon}\,\left|\left\langle  R_3\left(p,\tilde\epsilon \right) \,,\, n_{p} \right\rangle \right|  
%    \label{eq: Bounds on Gamma and R} \,.
%\end{multline}
To overcome the obstacle that both $\Gamma^{*}$ and $R^{*}$ are dependent on the point $p$, we appeal to compactness. Due to the Heine-Borel Theorem \cite[Theorem 3-40]{Apostol} we can cover the compact set $\mc M \cap \Omega$ by a finite number $k$ of closed balls $\mc B_{i}$, where $i \in \{1,\dots,k\}$. Since $f_{x}(p)$ and $F_{x}(p)$ are continuous, $\Gamma(p)$ is a continuous function of $p$. Thus on each of these balls $\Gamma^*$ and $R^*$ attain their minima and maxima as
$
\Gamma^*_{i} := \min_{p \in  \mc B_{i} \cap \mc M \cap \Omega} \Gamma^*
$
and
$
R^*_{i} := \max_{p \in \mc B_{i} \cap \mc M \cap \Omega} R^*
$, 
where $\Gamma_i^*$ and $R_i^*$ depend on $\mc B_{i} \cap \mc M \cap \Omega$. We define $\epsilon^*_{i} > 0$ such that the following inequality holds:
\begin{equation*} 
    \fao{p \in \mc B_{i} \cap \mc M \cap \Omega} \;
    \frac{1}{2}  \fvec c^T \Gamma^*_{i} \fvec c 
    >
    \epsilon^{*}_{i}\, R^*_i \,.
\end{equation*}
Therefore, we obtain together with \eqref{eq: Bounds on Gamma and R} that
\begin{equation*}
    \exo{\epsilon^*_{i} > 0}
    \fao{p \in \mc B_{i} \cap \mc M \cap \Omega} \;
    \frac{1}{2}  \fvec c^T \Gamma(p) \fvec c 
    > 
    \frac{1}{\epsilon^*_{i}} \left|\left\langle R_3\left(p,\epsilon^*_{i} \right) , n_{p}  \right\rangle \right|.
\end{equation*}
Because the number of balls is finite, we define $\epsilon^*>0$ as $\epsilon^* := \min_{i=1,\dots,k}\epsilon^*_{i}$ and have the result
\begin{equation*}
    \exo{\epsilon^* > 0}
    \fao{p \in \mc M \cap \Omega} \;
    \frac{1}{2} \fvec c^T  \Gamma(p) \fvec c 
    \;>\; 
    \frac{1}{\epsilon^*} \left|\left\langle R_3\left(p,\epsilon^* \right) \right\rangle , n_{p}  \right| 
    .
\end{equation*}
Thus $\epsilon^*$ provides a uniform bound for which the inner product \eqref{eq: Taylor linearization of inner product 2} is positive for every $ p \in \mc M \cap \Omega $. Clearly, the inner product is then also positive for every $p \in \mc M \cap \Omega$ if we choose any $\bar \epsilon \in (0,\epsilon^{*}]$. In other words, for any $\bar \epsilon \in (0,\epsilon^{*}]$, $\mc M_{\bar \epsilon}$ is overflowing invariant in $\Omega$.
\end{proof}

Theorem \ref{Theorem: Overflowing invariance} provides a checkable condition on the overflowing invariance of $\mc M$ within the compact set $\Omega$. Under further conditions on $\Omega$, hyperbolic instability of $\mc M \cap \Omega$ can be established.

\begin{cor}\label{Corollary: Instability under Overflowing Invariance}
Under the assumptions of Theorem \ref{Theorem: Overflowing invariance} and the additional assumption that $\Omega$ is an invariant strict superset of $\mc M \cap \Omega$, for any $\epsilon \in (0,\epsilon^{*}]$ the set $\Omega \setminus \bar{\mc M}_{\epsilon}$ is invariant.
\end{cor}

\begin{proof}
The set $\Omega$ can be partitioned by the non-empty sets $\bar{\mc M}_{\epsilon} \cap \Omega$ and $\Omega \setminus \bar{\mc M}_{\epsilon}$. Let us establish a correspondence of overflowing invariance and the flow of the vector field: since $\mc M_{\epsilon}$ is overflowing invariant in $\Omega$, we have for every $x_{0} \in\partial \mc M_{\epsilon} \cap \Omega$ and for all $t > 0$ that 
$
\Phi(t,x_{0}) \not \in 
\bar{\mc M}_{\epsilon} \cap \Omega
$. 
Therefore, a trajectory starting off $\bar{\mc M}_{\epsilon} \cap \Omega$ is bounded away from the partition $\bar{\mc M}_{\epsilon} \cap \Omega$. Invariance of $\Omega \setminus \bar{\mc M}_{\epsilon}$ follows then immediately from the invariance of $\Omega$.
\end{proof}

Corollary \ref{Corollary: Instability under Overflowing Invariance} allows a straightforward instability check of the set $\mc M \cap \Omega$, simply by analyzing the linearized vector field in \eqref{eq: overflowing invariance condition}. In the case that $\mc M$ is the origin and $\Omega$ is some nontrivial set containing $\mc M$, equation \eqref{eq: overflowing invariance condition} reduces to the equation obtained by Lyapunov's first method when using the identity as the Lyapunov matrix. Note that the results of this section can also be reversed, leading to asymptotic stability of a manifold \cite{DoerflerDiplomaThesis}. In the following section the geometric method will be applied to the link dynamics to show instability of the set $\mc W_{e} \setminus \mc E_{e}$.

% -----------------------------------------------------------------------------------------------------%
%      Global Stability Analysis of the Triangular Formation
% -----------------------------------------------------------------------------------------------------%

\section{Global Stability Analysis of the Target Formation}\label{Section: Global Stability Analysis of the Triangular Formation}

% -----------------------------------------------------------------------------------------------------%
\subsection{Equilibria and Invariant Sets of the Link Dynamics}\label{Subsection: Equilibria of the Link Dynamics}

The limit set of the link dynamics $\mc W_{e}$ can from \eqref{eq: dot V(e)} be parametrized as 
$
\mc W_{e} = \{e \in \mathrm{Im}\hat H:\, e_{1}\psi_{1}=e_{2}\psi_{2}=e_{3}\psi_{3}\}
$,
which is the set of equilibria of the link dynamics \eqref{eq: e-dynamics}. Clearly, $\mc W_{e}$ contains besides the target formation $\mc E_{e}$ also the set of collinear (non-rigid) equilibria. Let the set of collinear links be termed the \textit{line set} $\mc N_{e}$. By equation \eqref{eq: link space} the three links are linearly dependent, and $\mc N_{e}$ is naturally parameterized by two links and the planar $90^{\circ}$ rotation matrix $J$: 
\begin{equation*}
	\mc N_{e} 
	=
	\{ e\in \mathrm{Im}\hat H : \, e_{1}^{T}J\,e_{2} = \fvec 0 \} \,,
	\mbox{ where $J=$\footnotesize{
	$
	\left[ \begin{array}{rr}
		 0 & 1 \\  -1 & 0
	\end{array} \right]  
	$}} \,.	
\end{equation*}
Note that $\mc E_{e}$ and $\mc N_{e}$ are a positive distance apart, which follows directly from Theorem \ref{Theorem: Semiglobal Exp. Stability}. It can easily be checked that $\mc N_{e}$ is invariant with respect to the link dynamics, which implies that initially collinear robots remain collinear for all time \cite{Cao1,Cao2,Anderson} and formation control fails. %Furthermore, note that $e \in \mc W_{e}\setminus \mc E_{e}$ implies that at least one $\psi_{i} \neq \fvec 0$ which together with \eqref{eq: link space} means that the robots are collinear. Thus formation control fails in 
%$
%\mc W_{e}\setminus \mc E_{e} 
%= 
%\mc N_{e} \cap \mc W_{e}
%$. 

% -----------------------------------------------------------------------------------------------------%
\subsection{Instability of the Line Set}\label{Subsection: Overflowing Invariance and Instability of the Line Set}

Our goal is to show that trajectories of the link dynamics are bounded away from the line set $\mc N_{e}$. References \cite{Cao1,Cao2,Cao2008,Smith} carry out a Lyapunov approach and show that a function related to the point-to-set distance to the line set $\mc N_{e}$ is locally increasing (near the collinear equilibria $\mc N_{e} \cap \mc W_{e}$). Up to a multiplicative constant the chosen Lyapunov functions are equivalent to the oriented area of the triangle, which is $1/2\,e_{1}^{T} J e_{2}$. Obviously, these Lyapunov functions are problem-specific for the triangular formation and do not extend to other examples. By decomposing $\mc N_{e}$ into submanifolds and applying the results of the previous section, an analogous result is provable by purely algebraic calculations of equation \eqref{eq: overflowing invariance condition} and without guessing a Lyapunov function. 

First, we consider a subset of $\mc N_{e}$, the set of collocated robots defined by the zero set 
$
\mc X_{e}
=
\{ e \in \mathrm{Im} \hat H :\, e = \fvec 0 \}
$. Since $\mc X_{e}$ is the origin of $\mbb R^{6}$, it is an embedded submanifold of $\mbb R^{6}$ located in $\mathrm{Im}\hat H$. Its normal space $N_e \mc X_e$ can easily be parametrized as 
\begin{equation*}
	N_e \mc X_e 
	=
	\mb{column}\spano{ \begin{m} -I_2 & 0 & I_2 \\ I_2 & -I_2 & I_2 \\ 0 & I_2 & I_2 \end{m} } \,,
\end{equation*}
where the first four columns are within the link space and the last two are orthogonal to it. We now apply Theorem \ref{Theorem: Overflowing invariance} to show overflowing invariance of  $\mc X_{e,{\mc \epsilon}_{\mc X}}$, the tubular $\epsilon_{\mc X}$  neighbourhood of $\mc X_{e}$. Together with Corollary \ref{Corollary: Instability under Overflowing Invariance} this guarantees hyperbolic instability of $\mc X_{e}$.

\begin{lem}\label{Lemma: Overflowing Invariance of Xe}
Consider $e_{0} \in \mathrm{Im}\hat H$ such that $\mc X_{e} \cap \Omega(V(e_{0})) \neq \emptyset$. There exists $\epsilon_{\mc X}^{*} > 0$ such that for every $\epsilon_{\mc X} \in (0,\epsilon_{\mc X}^{*}]$ the set 
$
\Omega(V(e_{0})) \setminus \bar{\mc X}_{e,\epsilon_{\mc X}}
$ 
is invariant.
\end{lem}

\begin{proof}
We calculate the matrix $\Gamma_{\mc X_{e}}$ from equation \eqref{eq: overflowing invariance condition} for the invariant set $\mc X_{e}$. The Jacobian of the vector field \eqref{eq: e-dynamics} evaluated on $\mc X_e$ is obtained as $\hat H \, \mathrm{diag} ( -d_i^2\,I_2 )$, and the first four columns of $N_{e} \mc X_{e}$ provide a basis for the normal space of $\mc X_{e}$ within the link space. Thus we obtain
\begin{equation*}
	\Gamma_{\mc X_{e}}
    	= 
	 \begin{m} 
		\left(2\,d_1^2 + 4\,d_2^2\right)\,I_2 		&  \left(d_1^2 - 3\, d_2^2 - d_3^2\right)\,I_2 \\ 
		\left(d_1^2 - 3\, d_2^2 - d_3^2 \right)\,I_2 	&   \left(2\,d_2^2 + 4\,d_3^2\right)\,I_2 
	\end{m} 
	\,.
\end{equation*}
%whose first principal minor is positive. If the $d_{i}$ satisfy the triangle inequalities the second principal minor can be lower-bounded by
%\begin{multline*}
%	\det \Gamma_{\mc X_{e}}
%	=
%	10(d_{1}^{2}d_{2}^{2} +d_{1}^{2}d_{3}^{2} + d_{2}^{2}d_{3}^{2}) - d_{1}^{4} - d_{2}^{4} - d_{3}^{4}
%	\\
%	>
%	10(d_{1}^{2}d_{2}^{2} +d_{1}^{2}d_{3}^{2} + d_{2}^{2}d_{3}^{2}) - d_{1}^{2}(d_{2}^{2} + d_{3}^{2}) - d_{2}^{2}(d_{1}^{2} + d_{3}^{2}) - d_{3}^{2}(d_{1}^{2} + d_{2}^{2})
%	\,,
%\end{multline*}
%which is strictly positive. 
A simple argument shows that the principal minors of $\Gamma_{\mc X_{e}}$ are positive whenever $d_{1}$, $d_{2}$, and $d_{3}$ satisfy the triangle inequalities. Thus the assumptions of Theorem \ref{Theorem: Overflowing invariance} and Corollary \ref{Corollary: Instability under Overflowing Invariance} are satisfied within the compact and invariant set $\Omega(V(e_{0}))$, and the lemma follows immediately.
\end{proof}

In order to continue, consider the smooth function $F: \mbb R^{6} \to \mbb R^{3}$,
\begin{equation*}
	F(e) 
	= 
	\begin{m}  
		e_{1}^{T}\,J\,e_{2} \\
		e_{1} + e_{2} + e_{3} 
	\end{m} \,,
\end{equation*}
and note that $\mc N_{e}$ can be written as the zero set $\mc N_{e} = F^{-1}(\fvec 0)$. The Jacobian of $F(e)$  is given by
\begin{equation*}
	F_{e}(e)
	= 
	\begin{m}
		-e_2^T\,J  &  e_1^T \,J  & \fvec 0 \\
		I_2           & I_2            & I_2 
	\end{m}
\end{equation*}
and has constant rank three for all $e \in F^{-1}(\fvec 0) \setminus \{\fvec 0 \}$ and a rank loss for $e= \fvec 0$. Thus $\mc N_{e}$ is not a submanifold. However, if we subtract the set $\mc X_{e}$ together with the negatively invariant set $\mc X_{e,\epsilon_{\mc X}} \cup \partial \mc X_{e,\epsilon_{\mc X}}$, with $\epsilon_{\mc X}$ from Lemma \ref{Lemma: Overflowing Invariance of Xe}, then we obtain 
$
\mc N_{e}' := \mc N_{e} \setminus \bar{\mc X}_{e,\epsilon_{\mc X}}
$
as an embedded submanifold in $\mbb R^{6}$, which follows directly from the parameterization of $\mc N_{e}'$ via $F$ \cite[Proposition 5.28]{lee2003ism}.  Note that we have to be cautious in the later application of Theorem \ref{Theorem: Overflowing invariance} to $\mc N_{e}'$ since $\mc N_{e}'$ is neither open nor closed in the topology of $\mbb R^{6}$. Note also that $\mc N_{e}'$ is located in the link space, it is invariant, due to hyperbolic instability of $\bar{\mc X}_{e,\epsilon_{\mc X}}$, and its normal space is parametrized by $\mathrm{Im}F_{e}(e) |_{\mc N_{e}'}$ and is well defined. Similar to $N_{e} \mc X_{e}$ above, the normal space $N_e \mc N_e' $ can be split into components orthogonal and parallel to $(I_{2},I_{2},I_{2})$, the normal vector of the link space. We refer to page 140 of the thesis \cite{DoerflerDiplomaThesis} for the easy calculations leading to the parameterization 
\begin{equation*}
    	N_e \mc N_e' 
	= 
	\mb{column}\spano{ \left. \begin{m} -J\,e_2 & I_2 \\ -J\,e_3 & I_2 \\ -J\,e_1 & I_2 \end{m} \right|_{\mc N_e'} } \,.
\end{equation*}

The following lemma shows that no trajectory can approach the collinear equilibria via $\mc N_{e,\epsilon_{\mc N'}}'$.

\begin{lem}\label{Lemma: Overflowing Invariance of Ne'}
Consider $e_{0} \in \mathrm{Im}\hat H$ such that $\mc N_{e}' \cap \mc W_{e} \neq \emptyset$. There exists an $\epsilon_{\mc N'}^{*} > 0$, such that for every $\epsilon_{\mc N'} \in (0,\epsilon^{*}_{\mc N'}]$ the tubular $\epsilon_{\mc N'}$ neighbourhood $\mc N_{e,\epsilon_{\mc N'}}'$ is overflowing invariant in $\mc W_{e}$.
\end{lem}

Before we continue to the proof of Lemma \ref{Lemma: Overflowing Invariance of Ne'}, we state the following algebraic relationship:

\begin{lem}\label{Lemma: Ming's Lemma on potential functions}\cite[Lemma 6]{Cao1}
For any $e \in \mc N_{e}\cap \mc W_{e} $ we have that $\psi_{1}+\psi_{2}+\psi_{3} < 0$.
\end{lem}

Lemma \ref{Lemma: Ming's Lemma on potential functions} can be proved by considering all possible cases of collinear and collocated robots. With this algebraic relationship we can now move on to the proof of Lemma \ref{Lemma: Overflowing Invariance of Ne'}.

\begin{proof}[Proof of Lemma \ref{Lemma: Overflowing Invariance of Ne'}]
First, we verify that $\mc N_{e}' \cap \mc W_{e}$ is closed. From Lemma \ref{Lemma: Overflowing Invariance of Xe} we know that $\mc X_{e}$ is hyperbolically unstable and that on
$
{\mc X}_{e,\epsilon_{\mc X}} \cup \partial {\mc X}_{e,\epsilon_{\mc X}}
$
the vector field is pointing outward and is thus strictly non-zero. In short, $\mc X_{e}$ is an isolated subset of the collinear equilibria $\mc N_{e} \cap \mc W_{e}$. Due to continuity of the vector field, there can be no equilibrium set, such as $\mc N_{e}' \cap \mc W_{e}$, arbitrarily close to $\mc N_{e}'$'s boundary $\mc N_{e} \cap \partial {\mc X}_{e,\epsilon_{\mc X}}$. This proves that $\mc N_{e}' \cap \mc W_{e}$ is closed. Compactness follows from the fact that $\mc W_{e}$ is compact. The Jacobian of the vector field is given by $\hat H \, \mb{diag}( \Theta_{i})$ with 
$
\Theta_{i} :=
(\psi_i I_2 + 2 e_i e_i^T)|_{e \in \mc N_e' \cap \mc W_e }
$.
For notational convenience, the argument $e \in \mc N_e' \cap \mc W_e$ is left out in the following calculations. A basis for the normal space of $\mc N_{e}'$ within the link space is given by the first column of $N_{e} \mc N_{e}'$. Following an easy calculation we obtain term from \eqref{eq: overflowing invariance condition} as 
\begin{equation*}
    	\Gamma_{\mc N_{e}'}(e) 
	=
	e_2^T  J \Theta_1  J  e_2  -  e_2^T J \Theta_2 J e_3   +   e_3^T J \Theta_2 J e_3  
	- e_3^T J \Theta_3 J e_1   +  e_1^T J \Theta_3 J e_1 - e_1^T J \Theta_1 J e_2 \,.
\end{equation*}
%\begin{multline*}
%    	\Gamma_{\mc N_{e}'}(e) 
%	=
%	e_2^T  J \Theta_1  J  e_2  -  e_2^T J \Theta_2 J e_3   +   e_3^T J \Theta_2 J e_3  
%	\\
%	- e_3^T J \Theta_3 J e_1   +  e_1^T J \Theta_3 J e_1 - e_1^T J \Theta_1 J e_2 \,.
%\end{multline*}
The expression $J \Theta_i J$ simplifies further to
$
J \Theta_i J 
= 
-\psi_i I_2  + 2 J e_i e_{i}^{T} J 
$\,.
Note that for $e \in \mc N_e' \cap \mc W_e$ the links are collinear and thus we have for any $i,j,k \in \{1,2,3\}$ that
\begin{equation*}
	e_{j}^{T}\,J\,\Theta_i\,J\,e_{k} 
	= 
	-\psi_{i}\,e_{j}^{T}\,e_{k} + 2\,e_{j}^{T}\,J\,e_{i} \,e_{i}^{T}\,J\,e_{k} 
	=
	-\psi_{i}\,e_{j}^{T}\,e_{k} \,.
\end{equation*}
Therefore, $\Gamma_{\mc N_{e}'}(e)$ simplifies to
\begin{equation}
    \Gamma_{\mc N_{e}'}(e) 
    = 
    \left(  e_2^T \, \psi_1 \, e_1  + e_3^T \, \psi_2 \, e_2 + e_1^T \, \psi_3 \, e_3 \right) 
    \\
    - 
    \bigl( \psi_1 \, \vectornorm{e_2}^2  + \psi_2 \, \vectornorm{e_3}^2 + \psi_3 \, \vectornorm{e_1}^2 \bigr)  
    \label{eq: GammaN simplified} \,.
\end{equation}
%\begin{multline}
%    \Gamma_{\mc N_{e}'}(e) 
%    = 
%    \left(  e_2^T \, \psi_1 \, e_1  + e_3^T \, \psi_2 \, e_2 + e_1^T \, \psi_3 \, e_3 \right) 
%    \\
%    - 
%    \bigl( \psi_1 \, \vectornorm{e_2}^2  + \psi_2 \, \vectornorm{e_3}^2 + \psi_3 \, \vectornorm{e_1}^2 \bigr)  
%    \label{eq: GammaN simplified} \,.
%\end{multline}
Now we evaluate this expression on the compact set $\mc N_e' \cap \mc W_e$. Remember that for any $e \in \mc N_e' \cap \mc W_e$ it holds that 
$
e_1\,\psi_1 = e_2\,\psi_2 = e_3\,\psi_3
$.
Consequently, the first term of \eqref{eq: GammaN simplified} is zero:
\begin{equation*}
	e_2^T \, \psi_1 \, e_1  + e_3^T \, \psi_2 \, e_2 + e_1^T \, \psi_3 \, e_3 
	= 
	\left(e_1 + e_2 + e_3 \right)^{T} \, e_1\,\psi_1 
	= 
	0
	\,.
\end{equation*}
To analyze the second term we consider the cases where two or none of the robots are collocated:

\noindent\underline{case 1:}
    $e \in \left\{ \mc N_e' \cap \mc W_e \right\} \cap \{e \in \mathrm{Im} \hat H : \, e_{i}=\fvec 0 , e_{j \neq i} \neq \fvec 0 ,\,  i,j \in \{1,2,3\} \}$: 
    Suppose robot 1 and robot 2 are collocated, that is, $e_1 = \fvec 0$. It follows that $\psi_1 = -d_1^2 < 0$, $e_2=-e_3$ and also 
    $
    \fvec 0 = e_{1}\,\psi_{1} = e_2\,\psi_2 = e_3\,\psi_3 = -e_2\,\psi_3 = -e_{3} \,\psi_{2}
    $. 
    Thus we obtain from \eqref{eq: GammaN simplified} that
	\begin{equation*}
            \Gamma_{\mc N_{e}'}(e) 
            = 
            - \bigl( -d_1^2 \, \vectornorm{e_2}^2  + e_3^T \psi_2 \, e_3 +\, \fvec 0 \bigr) 
            \;\;=\;\; 
            d_1^2\,\vectornorm{e_2}^2 > 0\,.
    \end{equation*}
    The proof for $e_2 = \fvec 0$ and $e_3 = \fvec 0$ is analogous.

\noindent\underline{case 2:}
     $e \in \left\{ \mc N_e' \cap \mc W_e \right\} \cap \{e \in \mathrm{Im} \hat H : \, e_{i} \neq \fvec 0 \,,\, i \in \{1,2,3\} \}$: 
     Suppose all three robots are collinear but none of them are collocated. Then there exists $x \in \mbb R \setminus \{-1,0\}$ such that
     $
     e_2 = x \, e_1
     $ 
     and 
     $
     e_3 = -e_{1} - e_{2} = - (1+x) \, e_1  
     $\,.
     It follows then with $e_{1}\,\psi_{1} = e_2\,\psi_2 = e_3\,\psi_3$ that
	$
	\psi_2 
	= 
	\psi_{1}/x
	$ 
	and 
	$
	\psi_3 = - \psi_1/(1+x)
	$,
	and from Lemma \ref{Lemma: Ming's Lemma on potential functions} we get the condition $\psi_{1}\lambda(x)<0$, where 
	$
	\lambda(x) 
	:= 
	1 + 1/x - 1/(1+x)
	$.
    After some algebraic manipulations we can reformulate \eqref{eq: GammaN simplified} in terms of $e_{1}$, $\psi_{1}$, $x$, and $\lambda(x)$ as a product of strictly positive terms:
     \begin{equation*}
       		 \Gamma_{\mc N_{e}'}(e)
       		 = 
		 2  \vectornorm{e_1}^2 
		 \cdot 
		 ( - \psi_1 \,  \lambda^{-1}(x) )
		 \cdot
		 ( ( x + 1/2 )^2 + 3/4 )^3 / (x^2 (1+x)^2 ) 
		 > 0
		 \,.
       \end{equation*}
       
In summary, $\Gamma_{\mc N_{e}'}(e)>0$ for any $e$ in the compact set $\mc N_e' \cap \mc W_e$. Equivalently, there exists $\epsilon^{*}_{\mc N'} > 0$ such that for every $\epsilon_{\mc N'}  \in (0,\epsilon^{*}_{\mc N'}]$, $\mc N_{e,\epsilon_{\mc N'} }'$ is overflowing invariant in $\mc W_{e}$.
\end{proof}

From Lemma \ref{Lemma: Overflowing Invariance of Ne'} we conclude that the vector field \eqref{eq: e-dynamics} is pointing strictly outward on the set
\begin{equation*}
	\mc S_{\epsilon_{\mc N'}}
	:=
	\bigl\{ 
	\tilde e \in \mathrm{Im}\hat H:
	\tilde e 
	=
	e + \epsilon_{\mc N'} \cdot n_{e}  , e \in \mc N_{e}' \cap \mc W_{e} , 
	n_{e} \in N_{e} \mc N_{e}' , \norm{n_{e}} = 1 , \epsilon_{\mc N'}  \in (0,\epsilon^{*}_{\mc N'}]
	\bigr\}
	,
\end{equation*}
%\begin{multline*}
%	\mc S_{\epsilon_{\mc N'}}
%	:=
%	\bigl\{ 
%	\tilde e \in \mathrm{Im}\hat H:
%	\tilde e 
%	=
%	e + \epsilon_{\mc N'} \cdot n_{e}  , e \in \mc N_{e}' \cap \mc W_{e} , 
%	\\
%	n_{e} \in N_{e} \mc N_{e}' , \norm{n_{e}} = 1 , \epsilon_{\mc N'}  \in (0,\epsilon^{*}_{\mc N'}]
%	\bigr\}
%	,
%\end{multline*}
that is, the set of non-collinear links which can be reached from the equilibria $\mc N_{e}' \cap \mc W_{e}$ by going $\epsilon^{*}_{\mc N'}$ or less normally to $\mc N_{e}'$. After the simple but tedious algebraic calculations in the proofs of Lemma \ref{Lemma: Overflowing Invariance of Xe} and Lemma \ref{Lemma: Overflowing Invariance of Ne'}, we are now in a position to state our final result: 

\begin{theo}\label{Theorem: Overflowing invariance of Me}
Consider $e_{0} \in \mathrm{Im}\hat H$ such that $\mc N_{e} \cap \Omega(V(e_{0})) \neq \emptyset$. There exists an $\epsilon^{*} > 0$, such that for every $\epsilon \in (0,\epsilon^{*}]$ the set 
$
\Omega(V(e_{0})) 
\setminus 
\{
\mc N_{e}
\cup
\bar{\mc S}_{\epsilon}
\cup 
\bar{\mc X}_{e,\epsilon}
\}
$ is invariant.
\end{theo}

\begin{proof}
Let 
$
\epsilon^{*} = 
\min\{ \epsilon^{*}_{\mc X} , \epsilon^{*}_{\mc N'}\}
$ 
and let $\epsilon \in (0,\epsilon^{*}]$ be fixed. By Theorem \ref{Theorem: Semiglobal Exp. Stability}, for any initial condition $e_{0}$ the corresponding trajectory $\Phi(t,e_{0})$ is bounded in $\Omega(V(e_{0}))$ and will converge to a limit set in $\mc W_{e} = \mc E_{e} \cup \{ \mc N_{e} \cap \mc W_{e} \}$. Assume that trajectories starting off $\mc N_{e}$ approach the collinear equilibria $\mc N_{e} \cap \mc W_{e}$. These trajectories cannot first converge to $\mc N_{e} \setminus \mc W_{e}$ (in finite time) and then approach $\mc N_{e} \cap \mc W_{e}$ since then trajectories would intersect the invariant set $\mc N_{e}$ in non-equilibria. Furthermore, according to Lemma \ref{Lemma: Overflowing Invariance of Ne'}, a trajectory starting off 
$
\mc N_{e} 
\cup
\bar{\mc S}_{\epsilon}
$ 
cannot approach 
$
\mc N_{e}' \cap \mc W_{e}
$ 
via a neighbourhood of $\mc N_{e}'$ because it cannot enter $\bar{\mc S}_{\epsilon}$. By Lemma \ref{Lemma: Overflowing Invariance of Xe}, the set
$
\Omega(V(e_{0})) \setminus \bar{ \mc X}_{e,\epsilon}
$ 
is invariant, too. Therefore, a trajectory $\Phi(t,\xi_{0})$ with 
$
\xi_{0} \in  
\Omega(V(e_{0})) 
\setminus 
\{
\mc N_{e}
\cup
\bar{\mc S}_{\epsilon}
\cup 
\bar{\mc X}_{e,\epsilon}
\}
$ 
cannot approach 
$
\{ \mc N_{e} \setminus W_{e} \}
\cup
\{ \mc N_{e}'  \cap \mc W_{e} \} \cup \mc X_{e} 
=
\mc N_{e}
$. In particular, $\Phi(t,\xi_{0})$ will be bounded away from 
$
\mc N_{e}
\cup
\bar{\mc S}_{\epsilon}
\cup 
\bar{\mc X}_{e,\epsilon}
$. 
Finally note that $\epsilon$ can be chosen arbitrarily in $(0,\epsilon^{*}]$.
\end{proof}

Theorem \ref{Theorem: Overflowing invariance of Me} implies that initially not collinear robots will never be collinear and the corresponding trajectory will be bounded a strictly positive distance away from the collinear equilibria. By standard arguments \cite{Cao1,Cao2,Cao2008,DoerflerDiplomaThesis}, it can now be shown that a trajectory starting off $\mc N_{e}$ enters the level set $\Omega(\rho)$ from Theorem \ref{Theorem: Semiglobal Exp. Stability} within a finite time. Thus the target formation $\mc E_{e}$ is exponentially stable with $\mathrm{Im}\hat H \setminus \mc N_{e}$ as exact region of attraction. Spoken differently, initially not collinear robots converge exponentially to the specified triangular formation.

We conclude by discussing three possible extensions of the presented global stability analysis. 

(i) The final result in Theorem \ref{Theorem: Overflowing invariance of Me} can also be proved for more general and non-quadaratic potential functions, such as the  potential functions defined in \cite{Cao2} growing infinitely as two robots approach each other. Invariance of $\mathrm{Im}\hat H \setminus \mc X_{e}$ follows by standard Lyapunov arguments, Lemma \ref{Lemma: Ming's Lemma on potential functions} still holds \cite[Lemma 5]{Cao2}, and thus Lemma \ref{Lemma: Overflowing Invariance of Ne'} and Theorem \ref{Theorem: Overflowing invariance of Me} can be proven analogously. 

(ii) Switching topologies can be considered that are infinitesimally rigid, for example a cyclic topology with reverse link orientations, an undirected or acyclic topology. For each of these topologies local stability of $\mc E_{e}$ is guaranteed by Theorem \ref{Theorem: Semiglobal Exp. Stability} (as shown in \cite{DoerflerECC09}) with the exception of the acyclic topology which has to be analyzed based on its cascade structure \cite{Krick08,Cao2008}. Note that for each topology the same invariant set $\mc N_{e}$ arises and the manifold parameterizations are as before. However, for acyclic and undirected graphs the vector fields (and their Jacobians) are different and the positive definiteness of \eqref{eq: overflowing invariance condition} has to be verified separately for the two topologies. 

(iii) Higher order minimally rigid formations with undirected graphs have as limit sets also either the target formation $\mc E_{e}$ or non-rigid formations $\mc W_{e} \setminus \mc E_{e}$ \cite{DoerflerECC09}, which typically have collinear links. Therefore, the invariant sets are similar and can be analyzed with the methods presented.

\section{Conclusion}\label{Section: Conclusions}

The present paper considers a global stability analysis of the formation problem for autonomous robots. Based on the notion of overflowing invariance and geometric arguments a condition is derived in order to show instability of embedded submanifolds. This geometric method is then successfully applied to the example of a triangular formation with cyclic sensor graph in order to rule out undesired non-rigid limit sets of the closed-loop dynamics. The result relies on purely algebraic calculations and not on guessing a problem-specific Lyapunov function. 

%Although the presented geometric method is straightforward, the calculations can be tedious. However, it may be possible to relax our conditions by considering ellipsoidal neighbourhoods and render equation \eqref{eq: overflowing invariance condition} similar to the one obtained by  Lyapunov's first method. Our approach could then serve as a helpful tool in the global stability analysis of rigid formations.

%%%%%%%%%%%%%%%%%%%%%%%%%%%%%%%%%%%%%%%%%%%%%%%%%%%%%%%%%%%%%%%%%%%%%%%%%

% -----------------------------------------------------------------------------------------------------%
%      Bibliography
% -----------------------------------------------------------------------------------------------------%

% references are in file  ./references.bib
\bibliographystyle{IEEEtran}
\bibliography{references}

\begin{thebibliography}{10}
\providecommand{\url}[1]{#1}
\csname url@samestyle\endcsname
\providecommand{\newblock}{\relax}
\providecommand{\bibinfo}[2]{#2}
\providecommand{\BIBentrySTDinterwordspacing}{\spaceskip=0pt\relax}
\providecommand{\BIBentryALTinterwordstretchfactor}{4}
\providecommand{\BIBentryALTinterwordspacing}{\spaceskip=\fontdimen2\font plus
\BIBentryALTinterwordstretchfactor\fontdimen3\font minus
  \fontdimen4\font\relax}
\providecommand{\BIBforeignlanguage}[2]{{%
\expandafter\ifx\csname l@#1\endcsname\relax
\typeout{** WARNING: IEEEtran.bst: No hyphenation pattern has been}%
\typeout{** loaded for the language `#1'. Using the pattern for}%
\typeout{** the default language instead.}%
\else
\language=\csname l@#1\endcsname
\fi
#2}}
\providecommand{\BIBdecl}{\relax}
\BIBdecl

\bibitem{Eren2002}
T.~Eren, P.~Belhumeur, B.~Anderson, and S.~Morse, ``A framework for maintaining
  formations based on rigidity,'' \emph{{Proceedings of the 15th IFAC World
  Congress , Barcelona}}, 2002.

\bibitem{OlfatiSaber02}
R.~Olfati-Saber and R.~M. Murray, ``Distributed cooperative control of multiple
  vehicle formations using structural potential functions,'' \emph{{Proceedings
  of the 15th IFAC World Congress , Barcelona}}, 2002.

\bibitem{Hendrickx}
J.~Hendrickx, B.~Anderson, J.~Delvenne, and V.~Blondel, ``Directed graphs for
  the analysis of rigidity and persistence in autonomous agent systems,''
  \emph{International Journal of Robust and Nonlinear Control}, vol.~17,
  no.~10, pp. 960--981, 2007.

\bibitem{Anderson09}
B.~Anderson, C.~Yu, B.~Fidan, and J.~Hendrickx, ``Rigid graph control
  architectures for autonomous formations,'' \emph{IEEE Control Systems
  Magazine}, vol.~28, no.~6, pp. 48--63, 2008.

\bibitem{Krick08}
{L. Krick and M. Broucke and B. Francis}, ``{Stabilization of infinitesimally
  rigid formations of multi-robot networks},'' \emph{International Journal of
  Control}, pp. 423 -- 439, 2009.

\bibitem{DoerflerECC09}
F.~D{\"o}rfler and B.~Francis, ``{Formation Control of Autonomous Robots Based
  on Cooperative Behavior},'' \emph{{Proceedings of the 2009 European Control
  Conference, Budapest Hungary}}, 2009.

\bibitem{Anderson}
B.~Anderson, C.~Yu, S.~Dasgupta, and S.~Morse, ``Control of three co-leader
  formation in the plane,'' \emph{System and Control Letters}, pp. 573--578,
  2007.

\bibitem{Cao1}
M.~Cao, C.~Yu, S.~Morse, B.~Anderson, and S.~Dasgupta, ``Controlling a
  {T}riangular {F}ormation of {M}obile {A}utonomous {A}gents,'' in \emph{46th
  IEEE Conference on Decision and Control (CDC), New Orleans, LA, USA},
  December 2007.

\bibitem{Cao2}
------, ``Generalized {C}ontroller for {D}irected {T}riangle {F}ormations,'' in
  \emph{17th International Federation of Automatic Control World Congress
  (IFAC), Seoul, Korea}, July 2008.

\bibitem{Cao2008}
M.~Cao, B.~Anderson, S.~Morse, and C.~Yu, ``Control of acyclic formations of
  mobile autonomous agents,'' \emph{47th IEEE Conference on Decision and
  Control, Cancun, Mexico}, December 2008.

\bibitem{Smith}
S.~Smith, M.~Broucke, and B.~Francis, ``{Stabilizing a multi-agent system to an
  equilateral polygon formation},'' in \emph{Proceeds of the 17th International
  Symposium on Mathematical Theory of Networks and Systems (MTNS2006)}, 2006,
  pp. 2415--2424.

\bibitem{lee2003ism}
J.~Lee, \emph{{Introduction to Smooth Manifolds}}.\hskip 1em plus 0.5em minus
  0.4em\relax Springer, 2000.

\bibitem{Wiggins}
S.~Wiggins, \emph{Normally Hyperbolic Invariant Manifolds in Dynamical
  Systems}.\hskip 1em plus 0.5em minus 0.4em\relax Springer, 1991.

\bibitem{MultiVariableCalculus}
L.~Corwin, \emph{Multivariable Calculus}.\hskip 1em plus 0.5em minus
  0.4em\relax CRC Press, 1982.

\bibitem{Apostol}
T.~M. Apostol, \emph{Mathematical Analysis}.\hskip 1em plus 0.5em minus
  0.4em\relax Addison-Wesley Publishing Company Inc., 1957, 4th printing, 1964.

\bibitem{DoerflerDiplomaThesis}
F.~D{\"o}rfler, ``{G}eometric {A}nalysis of the {F}ormation {C}ontrol {P}roblem
  for {A}utonomous {R}obots,'' Diploma {T}hesis, University of Toronto, 2008.

\end{thebibliography}

% important: compile via LaTex -> LaTex -> BibTex -> LaTex !!!

%%%%%%%%%%%%%%%%%%%%%%%%%%%%%%%%%%%%%%%%%%%%%%%%%%%%%%%%%%%%%%%%%%%%%%%%%

\end{document}